\documentclass[12pt,leqno]{article}

\usepackage{amsfonts}
\usepackage{amssymb}

\newcommand{\cD}{{\cal D}}
\newcommand{\cE}{{\cal E}}

\newcommand{\cI}{{\cal I}}

\newcommand{\cO}{{\cal O}}

\newcommand{\ord}{{\rm ord}}

\newcommand{\Aut}{{\rm Aut}}

\newcommand{\Ind}{{\rm Ind}}
\newcommand{\Res}{{\rm Res}}

\newcommand{\Trace}{{\rm Trace}}
\newcommand{\Tr}{{\rm Tr}}
\newcommand{\w}{{\rm w}}
\newcommand{\tm}{{\rm t}}
\newcommand{\Cov}{{\rm Cov}}
\newcommand{\ram}{{\rm ram}}
\newcommand{\Gal}{{\rm Gal}}

\newcommand{\ZZ}{{\mathbb Z}}

\newcommand{\RR}{{\mathbb R}}
\newcommand{\NN}{{\mathbb N}}
\renewcommand{\AA}{{\mathbb A}}
\newcommand{\PP}{{\mathbb P}}
\newcommand{\QQ}{{\mathbb Q}}

\newcommand{\wm}{{\mathfrak m}}

\newcommand{\ra}{\rightarrow}
\def\rightepi{{\longrightarrow \kern-0.7em \rightarrow}}
\newcommand{\oplusm}{\mathop{\oplus}\limits}

\newcommand{\notteilt}{{\,\not{\kern-0.075em|}\,}}
\def\antiddots{\mathinner{\mkern1mu\raise1pt\vbox{\kern7pt\hbox{.}}\mkern2mu
    \raise4pt\hbox{.}\mkern2mu\raise7pt\hbox{.}\mkern1mu}}

\begin{document}

\vspace*{15ex}

\begin{center}
{\LARGE\bf Galois Structure of Zariski Cohomology \\
for Weakly Ramified Covers \\
\vspace*{0.2cm}
of Curves}\\

\bigskip
{\sc Bernhard K\"ock}
\end{center}

\bigskip

\begin{quote}
{\footnotesize {\bf Abstract}. We compute equivariant Euler
characteristics of locally free sheaves on curves, thereby
generalizing several results of Kani and Nakajima. For instance,
we extend Kani's computation of the Galois module structure of the
space of global meromorphic differentials which are logarithmic
along the ramification locus from the tamely ramified to the
weakly ramified case.



}

\end{quote}

\bigskip

\bigskip

\section*{Introduction}

Let $X$ be a smooth projective curve defined over an algebraically
closed field $k$ of characteristic $p$, and let $G \subseteq
\Aut(X/k)$ be a finite subgroup of automorphisms of $X$. The goal
of this paper is to compute the Galois module structure of the
Zariski cohomology groups of $X$ with values in an equivariant
locally free sheaf $\cE$ on $X$ such as the ideal sheaf of a
$G$-stable finite subset of points on $X$ or the sheaf of
differentials on $X$.

Our first result, see Theorem~(3.1), is an explicit formula for
the equivariant Euler characteristic
\[\chi(G,X,\cE) := [H^0(X,\cE)] - [H^1(X,\cE)]\]
considered as an element of the Grothendieck group $K_0(G,k)$ of
all finitely generated modules over the group ring $k[G]$. It
describes $\chi(G,X,\cE)$ in terms of the rank and degree of
$\cE$, the genus of the quotient curve $Y:=X/G$, the order of $G$
and of the higher ramification groups $G_{P,s}$, $P \in X$, $s \ge
0$, and the representations of the decomposition group $G_P$ on
the fibre $\cE(P)$ and on the cotangent space $\wm_P/\wm_P^2$ for
$P\in X$.

In the case the cover $\pi:X \ra Y$ is tamely ramified, this
formula becomes Theorem~1.1 in \cite{Ko2}, see Remark~(3.2). In
particular, if the order of $G$ is not divisible by $p$, it
implies the main result of the paper \cite{EL} by Ellingsrud and
L\o nsted, which in turn generalizes the classical Chevalley-Weil
formula, see \cite{Ko2}. While Theorem~(3.1) of the present paper
has the advantage of being available in general, i.e.\ without any
assumption on the ramification of $\pi$ or on the group $G$, it
has the disadvantage of computing the equivariant Euler
characteristic only in the ``weak'' Grothendieck group $K_0(G,k)$,
i.e.\ it yields only composition factors. In particular, if $p$
divides the order of $G$, we need a further input to describe the
actual $k[G]$-isomorphism class of the cohomology groups
$H^0(X,\cE)$ and $H^1(X,\cE)$, even if one of them vanishes.

Such an input is provided by Theorem~(1.1) parts of which may
already be found in the literature, see Remark~(1.2). It gives a
criterion for any fractional ideal in a local Galois extension to
have a normal basis element. In particular, it yields the
following fact, see the proof of Theorem~(2.1)(a): Let $\pi$ be
weakly ramified and let $D=\sum_{P\in X} n_P [P]$ be an
equivariant divisor on $X$ such that $n_P \equiv -1$ mod $e_P^\w$
for all $P\in X$ (where $e_P^\w$ denotes the order of the (first)
ramification group $G_{P,1}$); then the direct image
$\pi_*(\cO_X(D))$ of the associated equivariant invertible sheaf
$\cO_X(D)$ is a locally free $\cO_Y[G]$-module. (As in Erez' paper
\cite{Er}, the notion {\em weakly ramified} means that both tame
and the simplest kind of wild ramification are allowed, more
precisely, that all second ramification groups $G_{P,2}$, $P\in
X$, are trivial.) Using a standard argument in geometric Galois
module theory (see Chinburg's paper \cite{Ch} for the version most
suitable for our purposes), we obtain from this fact that the
equivariant Euler characteristic $\chi(G,X,\cO_X(D))$ lies in the
image of the (injective) Cartan homomorphism $K_0(k[G]) \ra
K_0(G,k)$ and, moreover, that $H^0(X,\cO_X(D))$ and
$H^1(X,\cO_X(D))$ are projective $k[G]$-modules, if one of them
vanishes, see Theorem~(2.1)(a).

This observation applied to the divisor $E:= \sum_{P\in X}(e_P^\w
-1)[P]$ together with the above-mentioned ``weak'' formula for
$\chi(G,X,\cE)$ applied to $\cE = \cO_X(E)$ allows us to construct
a certain canonical projective $k[G]$-module $N$ (depending only
on the action of $G$ on $X$) which embodies a global relation
between the local data $\wm_P/\wm_P^2$, $P\in X$, where
$\wm_P/\wm_P^2$ means the cotangent space of $X$ at $P$ together
with the obvious action of the decomposition group $G_P$, see
Theorem~(4.3). Moreover, using our ``weak'' formula again, we
express $\chi(G,X,\cO_X(D))$ as an integral linear combination of
classes of explicit {\em projective} $k[G]$-modules one of which
is $N$, see Theorem~(4.5). These theorems generalize results of
Kani and Nakajima from the tamely ramified to the weakly ramified
case, see Remarks~(4.4) and (4.6). Our approach to these theorems
(described above) generalizes the one used in \cite{Ko2} (a
special case of which may already be found in Borne's thesis
\cite{Bo}), but it is quite different from the ones used by Kani
and Nakajima.

Finally we give the following two applications of these theorems.
Firstly, we compute the $k[G]$-module structure of the first
cohomology group $H^1(X,\cI(S))$ of $X$ with values in the ideal
sheaf $\cI(S)$ of any $G$-stable finite subset $S$ of $X$ which
contains all wildly ramified points, see Corollary~(4.7). We refer
the reader to Pink's paper \cite{Pi} for the significance of this
ideal sheaf in his proof of a $p$-adic Grothendieck-Ogg-Shafaravic
formula. Secondly, if $S$ contains not only all wildly ramified
but all ramified points, we prove that the direct sum of the
$k[G]$-module $N$ with the space $H^0(X,\Omega_X(S))$ of global
meromorphic differentials, which are logarithmic along $S$, is a
free $k[G]$-module, see Corollary~(4.8). This result generalizes
Theorem~2 in Kani's paper \cite{Ka} again from the tamely ramified
to the weakly ramified case. In the tamely ramified case, Kani
furthermore deduces the $k[G]$-isomorphism class of the space
$H^0(X,\Omega_X)$ of all global {\em holomorphic} differentials
from this result, see Theorem~3 in \cite{Ka}. It would be
interesting to know whether this can also be done in the weakly
ramified case.

At this point we moreover mention that the conditions {\em $\pi$
is weakly ramified} and {\em $n_P \equiv -1$ mod $e_P^\w$ for all
$P \in X$} are not only sufficient, but also necessary for
$\pi_*(\cO_X(D))$ being locally free over $\cO_Y[G]$, see
Theorem~(1.1), and, if the degree of $D$ is sufficiently large,
also for $H^0(X,\cO_X(D))$ being a projective $k[G]$-module, see
Theorem~(2.1)(b). Without the assumption on the degree of $D$, it
might be true that these conditions are necessary for
$R\Gamma(X,\cO_X(D))$ being quasi-isomorphic to a perfect complex
of $k[G]$-modules (see Question~(2.6)).

The reader may also wish to consult the paper \cite{Vi} by
Vinatier for the current state of the art in the Galois module
theory of weakly ramified extensions of number fields.

{\bf Acknowledgments}. I would like to thank Niels Borne, Boas
Erez, Ian Leary, Richard Pink and Vic Snaith for helpful
discussions and for their encouraging interest. In particular,
Niels Borne has drawn my attention to Nakajima's paper \cite{Na3}
which is crucial for Theorem~(2.1)(b), Ian Leary has helped me
with the important Lemma~(4.2), and Vic Snaith has provided a
first approach to a central part of Theorem~(1.1).

\bigskip

\section*{\S1 The Normal Basis Theorem for Fractional Ideals
in Local Galois Extensions}

In this section we explicitly describe those fractional ideals in
any finite Galois extension of local fields for which the normal
basis theorem holds, i.e., which are free (of rank 1) over the
corresponding group ring.

Let $L/K$ be a finite Galois extension of local fields with Galois
group $G$ and (positive) residue class field characteristic $p$.
The corresponding extensions of discrete valuation rings, of
maximal ideals and of (perfect) residue class fields are denoted
by $\cO_L/\cO_K$, $\wm_L/\wm_K$ and $\lambda/\kappa$,
respectively. For any $s \ge -1$, let $G_s$ denote the $s$th
ramification group of the extension $L/K$. We recall (see Chapitre
IV in \cite{Se}): The ramification groups form a chain
\[ G=G_{-1} \supseteq G_0 \supseteq G_1 \supseteq G_2 \supseteq
\ldots\] of normal subgroups of $G$, $G_0 = {\rm ker}(G \ra {\rm
Gal}(\lambda/\kappa))$ is the inertia subgroup, $G_0/G_1$ is a
cyclic group of order prime to $p$, $G_s/G_{s+1}$ is an abelian
group of exponent $p$ for $s \ge 1$, and $G_s$ is the trivial
group for $s$ sufficiently big. The extension $L/K$ is called {\em
weakly ramified} ({\em tamely ramified, unramified}), iff $G_s=0$
for $s=2$ ($s=1$, $s=0$, respectively).

The following theorem generalizes the well-known result of
E.~Noether that $L/K$ is tamely ramified, if and only if $\cO_L$
is free (of rank 1) over the group ring $\cO_K[G]$. It implies in
particular that $L/K$ is weakly ramified, if and only if $\wm_L$
is free over $\cO_K[G]$. Furthermore it describes all fractional
ideals in $L$ which are $\cO_K[G]$-free. For example, if $L/K$ is
not weakly ramified, then there does not exist any $\cO_K[G]$-free
fractional ideal in $L$ at all.

{\bf (1.1) Theorem.} {\em Let $b\in \ZZ$. Then the fractional
ideal $\wm_L^b$ of $L$ is free over $\cO_K[G]$, if and only if
$L/K$ is weakly ramified and $b \equiv 1$~mod~$|G_1|$.}

{\em (1.2) Remark.} \\
(a) The only-if-part of Theorem~(1.1) follows from Theorem~3 and
the corollary of Proposition~2 in Ullom's paper \cite{Ul3}. The
tame case, i.e. $|G_1|=1$, and the case $b=1$, $G=G_1$ of the
if-part of Theorem~(1.1) are proved in Theorem~1 and Theorem~2 in
his paper. Unfortunately, he does not state the general case of
the if-part of Theorem~(1.1) which is essential for this paper.
Though it certainly can be proved with the methods he has
developed in his papers \cite{Ul1}, \cite{Ul2} and \cite{Ul3}, we
here give a coherent and self-contained proof of Theorem~(1.1) for
the reader's convenience.\\
(b) In the geometric case, Pink has given a ``global proof'' for
the fact that $L/K$ is weakly ramified, if and only if $\wm_L$ is
$\cO_K[G]$-free (see Corollary~3.6 in \cite{Pi}); to be precise,
the if-direction is proved there only under the additional
assumption $G=G_1$. Note that, in his terminology, {\em weakly
ramified} means {\em of type 2}.\\
(c) Let $\ord(G)$ be odd. Then Theorem~(1.1) also implies Erez'
theorem that $L/K$ is weakly ramified, if and only if the
so-called square root of the inverse different is $\cO_K[G]$-free
(see Theorem~1 in \cite{Er}).

We will use the following propositions and lemma in the proof of
Theorem~(1.1).

{\bf (1.3) Proposition.} {\em Let $I$ be any fractional ideal of
$L$. Then $I$ is free over $\cO_K[G]$, if and only if it is
projective over $\cO_K[G]$.}

{\em Proof.} This follows from a theorem of Swan (see
Corollary~6.4 on p.~567 in \cite{Sw}). \hfill $\square$

We recall that a $\ZZ[G]$-module $M$ is called {\em
cohomologically trivial}, iff the Tate cohomology groups
$\hat{H}^i(U,M)$, $i \in \ZZ$, vanish for all subgroups $U$ of
$G$.

{\bf (1.4) Proposition.} {\em Let $M$ be any $\cO_K[G]$-module.
Then $M$ is projective over $\cO_K[G]$, if and only if $M$ is
projective over $\cO_K$ and cohomologically trivial.}

{\em Proof.} If the ring of coefficients $\cO_K$ is replaced by
$\ZZ$, this is Th\'eor\`eme~7 on p.~151 in \cite{Se}. The same
proof applies, if the ring of coefficients is any Dedekind domain.
See also Proposition~4.1(a) on p.~457 in \cite{Ch}. \hfill
$\square$

As usual, we denote the multiplication with the norm element
$\sum_{\sigma \in G} [\sigma] \in \ZZ[G]$ by $\Tr_{L/K}$, the
different of $L/K$ by $\cD_{L/K}$, and the ramification index of
$L/K$ by $e:=e_{L/K}$. Furthermore, for any $r \in \RR$, the
standard notation $[r]$ means the greatest integer less than or
equal to $r$.

{\bf (1.5) Lemma.} {\em Let $b\in \ZZ$. Then we have:\\
(a) $(\wm_L^b)^G = \wm_K^a$, where $a = 1 + [\frac{b-1}{e}]$. \\
(b) $\Tr_{L/K}(\wm_L^b) = \wm_K^{a'}$, where $a' = [\frac{{\rm
ord}(\cD_{L/K}) +b}{e}]$.}

{\em Proof.} Let $a \in \ZZ$ and $a' \in \ZZ$ be defined by
$(\wm_L^b)^G = \wm_K^a$ and $\Tr_{L/K}(\wm_L^b) = \wm_K^{a'}$,
respectively. The obvious relation
\[\wm_L^{ae} = \wm_K^a\cO_L \subseteq \wm_L^b \subset
\wm^{a-1}_K\cO_L = \wm_L^{(a-1)e}\] implies $a \ge \frac{b}{e} >
a-1$, hence $a = 1 + [\frac{b-1}{e}]$. This proves assertion (a).
Furthermore we have: $\cD_{L/K}^{-1} = \{x \in L: \Tr_{L/K}(x
\cdot \cO_L) \subseteq \cO_K\}$. Hence:
\[\wm_L^{-ea'+b} = \wm_K^{-a'} \wm_L^b \subseteq \cD_{L/K}^{-1}
\subset \wm_K^{-(a'+1)} \wm_L^b = \wm_L^{-e(a'+1)+b}.\] Thus
$-ea'+b \ge -{\rm ord}(\cD_{L/K}) > -e(a'+1) +b$, i.e., $a' =
[\frac{{\rm ord}(\cD_{L/K}) +b}{e}]$, as was to be shown. \hfill
$\square$

{\em Proof of Theorem~(1.1).} We first proof the if-part. By
Propositions~(1.3) and (1.4), it suffices to show that the Tate
cohomology groups $\hat{H}^i(U,\wm_L^b)$, $i \in \ZZ$, vanish for
any subgroup $U$ of $G$. Since the extension $L/L^U$ is again
weakly ramified and since the first ramification group $U_1$ of
the extension $L/L^U$ is contained in $G_1$ (see Proposition~2 on
p.~70 in \cite{Se}), it suffices to consider the case $U=G$. The
usual spectral sequence argument yields the following fact for any
normal subgroup $N$ of $G$: If $\hat{H}^i(N,\wm_L^b) = 0$ and
$\hat{H}^i(G/N,(\wm_L^b)^N) = 0$ for all $i \in \ZZ$, then also
$\hat{H}^i(G,\wm_L^b) = 0$ for all $i \in \ZZ$. Applying this fact
to the filtration $G \supseteq G_0 \supseteq G_1 \supseteq \{1\}$
of $G$, we are reduced to show that $\hat{H}^i(G_1,\wm_L^b)=0$,
that $\hat{H}^i(G_0/G_1,I)=0$ for any fractional ideal $I$ in
$L^{G_1}$ and that $\hat{H}^i(G/G_0,J)=0$ for any fractional ideal
$J$ in $L^{G_0}$ (for all $i\in \ZZ$). Applying this fact again to
a filtration of the abelian group $G_1$ (of exponent $p$) whose
successive quotients are cyclic groups of order $p$, we see that
it finally suffices to show that $\hat{H}^i(G,\wm_L^b) =0$ for all
$i \in \ZZ$ and $b\in \ZZ$ with $b \equiv 1$ mod $|G_1|$ in the
following three cases (note that, for any subgroup $U$ of $G_1$ of
order $p$, we have $(\wm_L^b)^U = \wm_{L^U}^a$ with $a \equiv 1
\textrm{ mod }
|G_1/U|$ by Lemma~(1.5)(a)):\\
(i) $G_1$ is cyclic of order $p$ and $G=G_1$, i.e., $L/K$ is
totally wildly ramified of order~$p$.\\
(ii) $G_1 = \{1\}$ and $G = G_0$, i.e., $L/K$ is totally tamely
ramified.\\
(iii) $G_0 =1$, i.e., $L/K$ is unramified.\\
We first consider the case~(i). Since $G$ is cyclic, it suffices
to show $\hat{H}^i(G,\wm_L^b) = 0$ for $i=0$ and $i=1$. Using
Hilbert's formula for the order of the different (see
Proposition~4 on p.~72 in \cite{Se}) and the congruence $b\equiv 1
\textrm{ mod } p$, we obtain:
\[\left[\frac{{\rm ord}(\cD_{L/K}) + b}{p}\right] =
\left[\frac{2(p-1) +b}{p}\right] = 1 +
\left[\frac{b-1}{p}\right].\] Hence, by Lemma~(1.5), we have
$\Tr_{L/K}(\wm_L^b) = (\wm_L^b)^G$, i.e., $\hat{H}^0(G,\wm_L^b)
=0$. To show $\hat{H}^1(G,\wm_L^b)= 0$, it suffices to show that
the Herbrand difference
\[h(\wm_L^b) := {\rm
length}_{\cO_K}(\hat{H}^0(G,\wm_L^b)) - {\rm
length}_{\cO_K}(\hat{H}^1(G,\wm_L^b))\] vanishes. Since $\wm_L^b
\otimes_{\cO_K} K \cong L$ is $K[G]$-free of rank~1, we can find a
free $\cO_K[G]$-submodule $M$ of $\wm_L^b$ such that $\wm_L^b/M$
is of finite length. Now $h(\wm_L^b)=0$ follows from the
well-known (and easy) facts that the Herbrand difference vanishes
for any free $\cO_K[G]$-module and for any $\cO_K[G]$-module of
finite length, and that the Herbrand difference is additive on
short exact
sequences.\\
In the cases~(ii) and~(iii), we have to show that $\hat{H}^i(G,I)=
0$ (for all $i\in \ZZ$) for any arbitrary fractional $I$ of $L$.
Let $\pi_K$ be a prime element in $\cO_K$. Then we have
$\hat{H}^i(G,I/\pi_KI) = 0$ for all $i$. Indeed, in the case~(ii)
this is a consequence of Corollaire~1 on p.~138 in \cite{Se},
since $I/\pi_KI$ is annihilated by a power of $p$ and the order of
$G$ is prime to $p$; in the case~(iii) it is a consequence of the
classical normal basis theorem and of Proposition~(1.4), since
$I/\pi_KI$ is isomorphic to $\lambda$ as a $\kappa[G]$-module. The
long exact sequence associated with the short exact sequence of
$G$-modules
\[ 0 \ra I \,\, \stackrel{\pi_K}{\ra} \,\, I \ra I/\pi_KI \ra 0\] now
shows that $\pi_K\hat{H}^i(G,I) = \hat{H}^i(G,I)$ for all $i \in
\ZZ $. Finally, Nakayama's Lemma implies that $\hat{H}^i(G,I)=0$
for all $i \in \ZZ$. This completes the proof of the
if-direction in Theorem~(1.1). \\
We now prove the only-if-direction. So, let $b \in \ZZ$ such that
$\wm_L^b$ is $\cO_K[G]$-free. We first show that $L/K$ is weakly
ramified. Let $s$ be the greatest integer such that $G_s \not=
\{1\}$. We may assume that $s \ge 1$. Then there is a cyclic
subgroup $U$ of $G_s$ of order $p$. By Proposition~(1.4) we have
$\hat{H}^0(U,\wm_L^b)=0$, i.e., $\Tr_{L/L^U}(\wm_L^b) =
(\wm_L^b)^U$. Furthermore, we have ${\rm ord}(\cD_{L/L^U}) =
(s+1)(p-1)$ by Proposition~2 on p.~70 and Proposition~4 on p.~72
in \cite{Se}. Now Lemma~(1.5) applied to the extension $L/L^U$
implies that
\[s+1+\left[\frac{b-1-s}{p}\right] = \left[\frac{(s+1)(p-1) +b }{p}\right] =
1+\left[\frac{b-1}{p}\right].\] Hence $s$ must be equal to $1$,
i.e., $L/K$ is weakly ramified. It remains to show that $b \equiv
1$~mod~$|G_1|$. By Proposition~(1.4) we have
$\hat{H}^0(G_1,\wm_L^b) = 0$, i.e., $\Tr_{L/L^{G_1}}(\wm_L^b) =
(\wm_L^b)^{G_1}$. Furthermore, we have ${\rm ord}(\cD_{L/L^{G_1}})
= 2 \cdot (|G_1| -1)$ by Proposition~4 on p.~72 in \cite{Se}. Now,
Lemma~(1.5) applied to the extension $L/L^{G_1}$ implies:
\[2+\left[\frac{b-2}{|G_1|}\right] =
1 + \left[\frac{b-1}{|G_1|}\right].\] Hence $b \equiv 1$ mod
$|G_1|$, as desired. Thus, the proof of Theorem~(1.1) is complete.
\hfill $\square$

\bigskip

\section*{\S 2 Projectivity of Zariski Cohomology}

Let $k$ be an algebraically closed field of characteristic $p >
0$, $X$ a connected smooth projective curve over $k$ and $G$ a
finite subgroup of the automorphism group $\Aut(X/k)$ of
order~$n$.

In this section we give sufficient resp.\ necessary conditions
under which the Zariski cohomology groups of $X$ with values in an
equivariant invertible $\cO_X$-module are projective over the
group ring $k[G]$.

Let $\pi: X \ra Y:= X/G$ denote the canonical projection, and let
$g_X$ resp.\ $g_Y$ denote the genus of $X$ resp.\ $Y$.
Furthermore, for any $P \in X$, the decomposition group $\{\sigma
\in G: \sigma(P)=P\}$ is denoted by $G_P$, the ramification index
of $\pi$ at the place $P$ by $e_P$, the higher ramification groups
(see Chapitre IV in \cite{Se}) by $G_{P,s}$, $s\ge 0$, the wild
part of the ramification index, i.e.~$|G_{P,1}|$, by $e_P^\w$ and
the tame part of the ramification index, i.e.~$|G_P/G_{P,1}|$, by
$e_P^\tm$. We say that $\pi$ is {\em weakly ramified}, iff
$G_{P,s}$ is trivial for $s \ge 2$ and all $P \in X$.

We denote the Grothendieck group of all finitely generated
$k[G]$-modules (resp., of all finitely generated projective
$k[G]$-modules) by $K_0(G,k)$ (resp., by $K_0(k[G])$). We recall
from classical representation theory that the set of isomorphism
classes of irreducible $k[G]$-modules (resp., of indecomposable
projective $k[G]$-modules) forms a basis of $K_0(G,k)$ (resp., of
$K_0(k[G])$) and that the Cartan homomorphism $K_0(k[G]) \ra
K_0(G,k)$ is injective.

We recall that a {\em locally free $G$-module on $X$} is a locally
free $\cO_X$-module $\cE$ together with $\cO_X$-isomorphisms
$g^*(\cE)\ra \cE$, $g\in G$, which satisfy the usual composition
rules. For instance, if $D= \sum_{P \in X} n_P [P]$ is an {\em
equivariant divisor} on $X$ (i.e., $n_{\sigma(P)}=n_P$ for all
$\sigma \in G$ and $P \in X$), then the $\cO_X$-module $\cO_X(D)$
is a locally free $G$-module on $X$ of rank $1$. The Zariski
cohomology groups $H^i(X, \cE)$, $i\ge 0$, are then
$k$-representations of $G$ in the obvious way. Let
\[\chi(G,X,\cE) := [H^0(X,\cE)] - [H^1(X,\cE)] \in K_0(G,k)\]
denote the {\em equivariant Euler characteristic of $X$ with
values in $\cE$}.

{\bf (2.1) Theorem.} {\em Let $D= \sum_{P \in X} n_P [P]$ be an
equivariant divisor on $X$. \\
(a) Let $\pi$ be weakly ramified and $n_P \equiv -1$ mod $e_P^\w$
for all $P \in X$; then there exists a bounded complex $L^*$ of
finitely generated projective $k[G]$-modules  such that the
$k[G]$-module $H^i(X,\cO_X(D))$ is isomorphic to the $i$th
cohomology module $H^i(L^*)$ for all $i \in \ZZ$; in particular
we have:\\
(i) The equivariant Euler characteristic $\chi(G,X,\cO_X(D)) \in
K_0(G,k)$ lies in the image of the Cartan homomorphism
$K_0(k[G]) \ra K_0(G,k)$. \\
(ii) If one of the two cohomology groups $H^i(X,\cO_X(D))$, $i=0,
1$, vanishes, then the other one is a projective
$k[G]$-module.\\
(b) Let $\deg(D) > 2g_X -2$. If the $k[G]$-module
$H^0(X,\cO_X(D))$ is projective, then $\pi$ is weakly ramified and
$n_P \equiv -1$ mod $e_P^\w$ for all $P \in X$. }

{\em Proof.} \\
(a) Theorem~(1.1) implies that, for any $P\in X$, the
$\cO_{Y,\pi(P)}[G_P]$-module $\cO_X(D)_P = \wm_P^{-n_P}$ is free
after completion. From Corollary~(76.9) on p.~533 in \cite{CR} we
obtain that this is even true without completion. Hence the direct
image $\pi_*(\cO_X(D))$ is a locally free $\cO_Y[G]$-module.
Furthermore we have $H^i(X,\cO_X(D)) = H^i(Y,\pi_*(\cO_X(D)))$ for
all $i \in \ZZ$. Now Theorem~1.1 on p.~447 and Proposition~4.1(a)
on p.~457 in Chinburg's paper \cite{Ch} imply the first assertion
of Theorem~(2.1)(a). Statement~(i) is an immediate consequence of
this assertion, and statement~(ii) can be derived from it as in
the proof of Theorem~2
in Nakajima's paper~\cite{Na1}.\\
(b) We first prove this in the case that $G$ is cyclic of order
$p$. We fix a point $P \in X$. Let $N$ denote the greatest integer
such that $G_{P,N}$ is not trivial. We may assume that $N \ge 1$.
By definition of the higher ramification groups we have $N+1 =
\ord_P(\sigma(x)-x)$ where $\sigma \in G \backslash \{1\}$ and $x$
is any prime element of the local ring $\cO_{X,P}$. From Theorem~1
on p.~86 in \cite{Na3} we obtain that the non-negative integers
\[m_j := \frac{N}{p} + \langle\frac{n_P-jN}{p}\rangle -
\langle\frac{n_P-(j-1)N}{p}\rangle, \quad j=1, \ldots, p-1, \] are
zero; here, for a rational number $a$, $\langle a \rangle$ denotes
the fractional part of $a$, i.e., $0 \le \langle a \rangle < 1$
and $a- \langle a \rangle$ is an integer. Since $N$ is not
divisible by $p$ (see Lemma~1 on p.~87 in \cite{Na3}), there is a
solution $j_0 \in \{0, \ldots, p-1\}$ of the congruence $n_P
\equiv (j_0-1)N$  mod $p$. If $j_0 \not= 0$, then one of the
integers $m_j$, $j=1, \ldots, p-1$, would be positive, namely
$m_{j_0} = \frac{N}{p} + \langle \frac{n_P-j_0N}{p}\rangle$. Hence
$j_0 =0$, i.e., $N \equiv -n_P$ mod $p$; therefore
\[m_j = \frac{N}{p} + \langle\frac{-jN}{p}-\frac{N}{p}\rangle -
\langle \frac{-jN}{p} \rangle \quad \textrm{for all } j=1, \ldots,
p-1.\] Let now $j_1 \in \{1, \ldots, p-1\}$ be a solution of the
congruence $-j_1N \equiv 1$ mod $p$. If $N \not\equiv 1$ mod $p$,
then
\[m_{j_1} = \frac{N}{p} + \langle \frac{1}{p} - \frac{N}{p}
\rangle - \frac{1}{p} = \frac{N}{p} +1 -
\langle\frac{N}{p}\rangle\] would be positive. Hence $N\equiv 1$
mod $p$; therefore
\[m_j = \frac{N}{p} + \langle \frac{-j-1}{p} \rangle - \langle
\frac{-j}{p}\rangle = \frac{N}{p} + \frac{p-j-1}{p} -
\frac{p-j}{p} = \frac{N-1}{p}\] for all $j=1, \ldots, p-1$. Hence
we have $N=1$ and $n_P \equiv -1$ mod $p$, as was to be shown.\\
We now consider the general case. We fix a point $P \in X$ with
$e_P^\w \not=1$. (If such a point does not exist, we are done.)
Let $N$ denote the greatest integer such that $G_{P,N}$ is not
trivial. We choose a cyclic subgroup $H$ of $G_{P,N}$ of
order~$p$. Let $\eta: X \ra Z:= X/H$ denote the corresponding
cover. Since $H^0(X,\cO_X(D))$ is also projective as a
$k[H]$-module and since, for any $s\ge 0$, the intersection
$G_{P,s} \cap H$ is the $s$th ramification group of the cover
$\eta$ at $P$ (see Proposition~2 on p.~70 in \cite{Se}), we obtain
from the case considered above that $N=1$ and that $n_{P'} \equiv
-1$ mod $f_{P'}$ for {\em all} $P' \in X$ (where $f_{P'}$ denotes
the ramification index of $\eta$ at $P'$). It remains to show that
$n_P \equiv -1$ mod $e_P^\w$. We prove this by induction on $n:=
\ord(G)$. The case $n=1$ is trivial. Let $g_Z$ denote the genus of
$Z$, and let the divisor $E$ on $Z$ be defined by the equality
$\cO_Z(E)= \eta_*^H(\cO_X(D))$ of subsheaves of the constant sheaf
$K(Z)$, the function field of $Z$. Then by Lemma~(1.5)(a), the
multiplicity of $E$ at any point $Q \in Y$ is
$-\left(1+\left[\frac{-n_{\tilde{Q}}-1}{f_{\tilde{Q}}}\right]\right)
= -1 + \frac{n_{\tilde{Q}}+1}{f_{\tilde{Q}}}$ where $\tilde{Q} \in
Y$ is any point in the fibre $\eta^{-1}(Q)$. Thus we have:
\begin{eqnarray*}
\lefteqn{p \cdot \deg(E) = p \cdot \sum_{Q \in Z}
\left(-1+\frac{n_{\tilde{Q}}+1}{f_{\tilde{Q}}}\right)=}\\
&=& \sum_{P' \in X} (-f_{P'} + n_{P'} +1)\\
&=& \deg(D) - \sum_{P' \in X}(f_{P'}-1)\\
&\ge& 2(g_X-1) - \sum_{P' \in X} 2(f_{P'} -1)\\
&=& p \cdot (2g_Z-2)
\end{eqnarray*}
by Hurwitz' formula (see Corollary~2 on p.~301 in \cite{Ha}) and
Hilbert's formula (see Proposition~4 on p.~72 in \cite{Se}).
Furthermore $H^0(Z,\cO_Z(E)) = H^0(X, \cO_X(D))^H$ is projective
as a $k[G_{P,1}/H]$-module. Now the induction hypothesis applied
to the Galois cover $Z \ra X/G_{P,1}$ (with Galois group
$G_{P,1}/H$) and to the divisor $E$ on $Z$ implies that $-1 +
\frac{n_P+1}{p} \equiv -1$ mod $\frac{e_P^\w}{p}$, hence $n_P
\equiv -1$ mod $e_P^\w$, as was to be shown. This completes the
proof of Theorem~(2.1). \hfill $\square$

{\bf (2.2) Corollary.} {\em The following three assertions are
equivalent:\\
(a) The cover $\pi$ is tamely ramified.\\
(b) For every equivariant divisor $D$ on $X$ satisfying
$H^1(X,\cO_X(D)) =0$, the $k[G]$-module $H^0(X,\cO_X(D))$ is
projective. \\
(c) There exists a divisor $D$ on $X$ of the form $D=\pi^*(E)$,
$E$ any divisor on $Y$, with $\deg(D) > 2g_X-2$ such that the
$k[G]$-module $H^0(X,\cO_X(D))$ is projective.}

{\em Proof.} This immediately follows from Theorem~(2.1), since,
for any divisor $D= \sum_{P\in X} n_P [P]$ on $X$ of the form
$D=\pi^*(E)$, we have $n_P \equiv 0$ mod $e_P$ for all $P \in X$.
\hfill $\square$

{\bf (2.3) Corollary.} {\em Let $S$ be any non-empty $G$-stable
set of closed points on $X$ containing all ramified points, and
let $\Omega_X(S) := \Omega_X \otimes \cO_X(\sum_{P\in S}[P])$
denote the sheaf of meromorphic differentials on $X$ which are
logarithmic along $S$. Then, $\pi$ is weakly ramified, if and only
if the $k[G]$-module $H^0(X,\Omega_X(S))$ is projective.}

{\em Proof.} By Th\'eor\`eme~2.34 on p.~44 in \cite{Bo}, there is
an equivariant divisor $K_X = \sum_{P \in X} n'_P [P]$ on $X$ such
that $\cO_X(K_X) \cong \Omega_X$. Then the divisor $D= K_X +
\sum_{P \in S} [P]$ satisfies $\deg(D) > 2g_X-2$ (by Example~1.3.3
on p.~296 in \cite{Ha}). Hence, Theorem~(2.1)(b) implies the
if-part of Corollary~(2.3). We now prove the only-if-part. We have
a natural short exact sequence
\[0 \ra \pi^*(\Omega_Y) \ra \Omega_X \ra \Omega_{X/Y} \ra 0\]
of coherent $G$-modules on $X$. Let the divisor $K_Y = \sum_{Q\in
Y} r_Q [Q]$ on $Y$ be defined by the equality $\pi^*(\Omega_Y) =
\cO_X(\pi^*(K_Y))$ of subsheaves of $\cO_X(K_X)$. If $\pi$ is
weakly ramified, Hilbert's formula (see Proposition~4 on p.~72 in
\cite{Se}) then implies that $n'_P -e_Pr_{\pi(P)}=
(e_P-1)+(e_P^\w-1)$ for all $P \in X$. In particular, the divisor
$D$ satisfies the condition of Theorem~(2.1)(a), and
statement~(ii) of Theorem~(2.1)(a) finally implies the
only-if-part of Corollary~(2.3). \hfill $\square$

{\em (2.4) Remark.} \\
(a) Corollary~(2.2) is Theorem~2 in \cite{Na3}.\\
(b) Let $S$ and $\Omega_X(S)$ be as in Corollary~(2.3). Kani
states in Theorem~1 of \cite{Ka} that $H^0(X,\Omega_X(S))$ is a
projective $k[G]$-module, if and only if $\pi$ is tamely ramified.
Since there exist weakly ramified covers which are not tamely
ramified (see the case $r=1$ in Example~(2.5) below), the
only-if-direction of Corollary~(2.3) contradicts the
only-if-direction of Kani's result which seems to be wrong. (More
precisely, the final displayed equation and the computation of
$m_y$ for $y \in Y_\ram$ in the proof of the only-if-part of that
theorem seem to be wrong.)

The following example (taken from Hasse's paper \cite{Has}) shows
that every theoretically possible ``type'' of ramification for
Galois covers of smooth projective curves can occur even under the
rather restrictive conditions that the base curve is the
projective line, that the Galois group is cyclic of order~$p$ and
that only one point is ramified.

{\em (2.5) Example.} Let $r\in \NN$ such that $p$ does not divide
$r$. Let $k(x,y)$ be the cyclic field extension of the rational
function field $k(x)$ of degree $p$ given by the Artin-Schreier
equation $y^p-y=x^r$. Let $\pi: C \ra \PP^1_k$ denote the
corresponding cover of nonsingular curves. Then $\pi$ is
unramified precisely over $\AA^1_k \subset \PP^1_k$, and, at the
unique point $P \in C$ lying over $\infty \in \PP^1_k$, the
greatest integer $N$ such that $G_{P,N}$ is not trivial is equal
to $r$. Furthermore, the genus of $C$ is equal to
$\frac{(r-1)(p-1)}{2}$.

{\em Proof (extracted from \cite{Has}).} The Galois group $G=
\Gal(k(x,y)/k(x))$ is generated by the automorphism $\sigma$ given
by $\sigma(y) = y+1$. Thus, clearly, $\pi$ is unramified over
$\AA^1_k$. Let $P \in C$ be any point with $\pi(P) = \infty$.
Since $x^{-1}$ is a prime element at $\infty \in \PP^1_k$, the
equality $y^p-y=x^r$ implies that $\ord_P(y) <0$ and, more
precisely, that
\[\ord_P(y) = \frac{\ord_P(y^p)}{p} = \frac{\ord_P(y^p-y)}{p} =
\frac{\ord_P(x^r)}{p} = \frac{-r \cdot \ord_P(x^{-1})}{p}.\]
Hence, the cover $\pi$ is ramified at $P$ and we have $\ord_P(y)
=-r$. We choose $u,v \in \NN$ such that $-ur + vp =1$. Then
$y^ux^{-v}$ is obviously a prime element at $P$ and we have:
\begin{eqnarray*}
\lefteqn{N+1 = \ord_P(\sigma(y^ux^{-v}) - y^ux^{-v}) =}\\& =&
\ord_P\left(\left((y+1)^u - y^u\right)x^{-v}\right)=
\ord_P(uy^{u-1}x^{-v}) = r+1,
\end{eqnarray*}
as was to be shown. Finally, the genus $g_C$ of $C$ is determined
by the formula $2(g_C-1) = -2p + (N+1)(p-1)$ which is a
consequence of Hurwitz' formula (see Corollary~2.4 on p.~301 in
\cite{Ha}) and Hilbert's formula (see Proposition~4 on p.~72 in
\cite{Se}). Hence $g_C=\frac{(r-1)(p-1)}{2}$, as was to be shown.
\hfill $\square$

{\em (2.6) Question.} Recall that a bounded complex of finitely
generated projective $k[G]$-modules is called a {\em perfect
complex}. The proof of Theorem~(2.1)(a) actually yields the
following slightly stronger result: If $\pi$ is weakly ramified
and $n \equiv -1$ mod $e_P^\w$ for all $P\ \in X$, then
$R\Gamma(X,\cO_X(D))$ is quasi-isomorphic to a perfect complex. If
$\deg(D) > 2g_X -2$, then Theorem~(2.1)(b) shows that also the
converse statement is true. Is the converse statement true in
general?

\bigskip

\section*{\S 3 Computing Equivariant Euler Characteristics}

In the paper \cite{Ko2} we have given a formula for the Euler
characteristic of equivariant Zariski sheaves on curves in the
tamely ramified case. In this section we generalize that formula
to the general case.

We keep the assumptions and notations introduced in Section~2. We
do not assume any condition on the ramification of the cover $\pi:
X \ra Y:=X/G$ in this section. The main focus in this paper is on
the case $p = {\rm char}(k) > 0$, but everything in this and the
following section is true and interesting (though classical) also
in the case $p=0$ (if the condition {\em prime to $p$} is regarded
as the empty condition and the term {\em $p$-group} means {\em
trivial group}).

Let $\cE$ be a locally free $G$-module on $X$ of rank $r$. For any
$P\in X$, we view the fibre $\cE(P):=\cE_P/{\mathfrak m}_P\cE_P$
as a $k$-representation of the decomposition group $G_P$.
Furthermore, the obvious representation of $G_P$ on the cotangent
space $\wm_P/\wm_P^2$ (or the corresponding character $G_P \ra
k^\times$) is denoted by $\chi_P$.

The following theorem computes the equivariant Euler
characteristic $\chi(G,X,\cE)$ of $X$ with values in $\cE$.

{\bf (3.1) Theorem.} {\em We have in $K_0(G,k)_\QQ$:
\begin{eqnarray*}
\lefteqn{\chi(G,X,\cE) =} \\
&=& \left((1-g_Y)r + \frac{1}{n} \deg(\cE) -
\frac{r}{2n}\sum_{P\in X}\left((e_P^\w-1)(e_P^\tm+1) +\sum_{s\ge
2}(|G_{P,s}|-1)\right)\right)[k[G]]\\
&& \mbox{} - \frac{1}{n} \sum_{P\in X} e_P^\w
\sum_{d=1}^{e_P^\tm-1} d \cdot \left[\Ind_{G_P}^G\left(\cE(P)
\otimes \chi_P^d\right)\right].
\end{eqnarray*}}

{\em (3.2) Remark.}\\
(a) If $\pi$ is tamely ramified, then the upper sum over $P\in X$
obviously vanishes. In particular, Theorem~(3.1) generalizes
Theorem~1.1 in \cite{Ko2}. In the next section, we will apply
Theorem~(3.1) in the weakly ramified case, i.e., when
$\sum_{s\ge 2}(|G_{P,s}|-1)$ vanishes for all $P \in X$.\\
(b) The reader may wish to compare Theorem~(3.1) with
Th\'eor\`eme~3.18 in Borne's paper \cite{Bo2} which gives an
alternative expression of the left hand side.

{\em Proof of Theorem~(3.1).} By classical representation theory
(see \S 82 in \cite{CR}) it suffices to show that both sides
coincide after restricting to any cyclic subgroup $H$ of $G$ of
order prime to $p$. Therefore it suffices to show that, by
restricting to $H$, we obtain the corresponding formula for the
action of $H$ on $X$ and that the character values (of the
original formula) coincide at any
element of $G$ whose order is prime to $p$. \\
We first show that character values at $\sigma =1$, i.e.\ the
$k$-dimensions, coincide. The $k$-dimension of the right hand side
is
\begin{eqnarray*}
\lefteqn{n(1-g_Y)r + \deg(\cE) - \frac{r}{2}\sum_{P\in
X}\left((e_P^\w-1)(e_P^\tm+1) + \sum_{s\ge 2}(|G_{P,s}|-1)\right)}
\hspace*{20em}\\
&& \mbox{}-\frac{1}{n} \sum_{P\in X}e_P^\w \sum_{d=1}^{e_P^\tm-1}
d \cdot r \cdot \frac{n}{e_P}.
\end{eqnarray*}
By Hurwitz' formula (see Corollary~2.4 on p.~301 in \cite{Ha}) and
Hilbert's formula for the order of the different (see
Proposition~4 on p.~72 in \cite{Se}), this is equal to
\begin{eqnarray*}
\lefteqn{(1-g_X)r + \frac{r}{2}\sum_{P\in X} \sum_{s\ge 0}
(|G_{P,s}|-1) +\deg(\cE)} \\
&& \mbox{} - \frac{r}{2} \sum_{P\in X} \left((e_P^\w
e_P^\tm - e_P^\tm +e_P^\w - 1) + \sum_{s\ge 2} (|G_{P,s}|-1)\right)\\
&& \mbox{} - \sum_{P\in X} \frac{r}{e_P^\tm} \cdot
\frac{e_P^\tm(e_P^\tm-1)}{2} =\\
&=& (1-g_X)r + \deg(\cE) + \frac{r}{2} \sum_{P\in X}(e_P^\tm -1)\\
&& \mbox{} - \frac{r}{2} \sum_{P\in X} (e_P^\tm-1)\\
&=& (1-g_X)r + \deg(\cE).
\end{eqnarray*}
This is equal to the $k$-dimension of the left hand side by the
theorem of Riemann-Roch (see \S 1 in Chapter~IV of \cite{Ha} and
Exercise~6.11 on p.~149 in \cite{Ha}).\\
We now fix an element $\sigma \in G\backslash\{1\}$ whose order is
prime to $p$ and show that the character values of both sides
coincide at $\sigma$. Using $\Trace(\sigma|k[G])=0$ and, more
generally, the well-known formula for an induced character (see
formula~(38.3) on p.~266 in \cite{CR}), we obtain for the
character value of the right hand side at $\sigma$:
\begin{eqnarray*}
\lefteqn{-\frac{1}{n}\sum_{P\in X} e_P^\w \sum_{d=1}^{e_P^\tm-1} d
\cdot \Trace\left(\sigma|\Ind_{G_P}^G\left(\cE(P)\otimes \chi_P^d\right)\right) =}\\
&=& -\frac{1}{n} \sum_{P\in X} \sum_{d=1}^{e_P^\tm-1} \frac{e_P^\w
\cdot d}{e_P} \sum_{\tau \in G,\,\, \tau^{-1}\sigma\tau \in G_P}
\Trace\left(\tau^{-1}\sigma\tau|\cE(P)\right) \cdot
\chi_P^d\left(\tau^{-1}\sigma\tau\right)
\\
&=& -\frac{1}{n} \sum_{P\in X}
\sum_{d=1}^{e_P^\tm-1}\frac{d}{e_P^\tm} \sum_{\tau \in G,\,\,
\sigma \in G_{\tau(P)}} \Trace\left(\sigma|\cE(\tau(P))\right)
\cdot
\chi_{\tau(P)}^d(\sigma) \\
&=& -\sum_{P\in X^\sigma} \Trace\left(\sigma|\cE(P)\right) \cdot
\frac{1}{e_P^\tm} \sum_{d=1}^{e_P^\tm-1}d \cdot \chi_P^d(\sigma)
\end{eqnarray*}
where $X^\sigma := \{P\in X: \sigma(P)=P\}$. Since $G_{P,1}$ is a
$p$-group (see Corollaire~3 on p.~75 in \cite{Se}), the character
$\chi_P: G_P \ra k^\times$ factors modulo $G_{P,1}$, and the
induced character $\bar{\chi}_P: G_P/G_{P,1} \ra k^\times$ is
injective by Corollaire~1 on p.~75 in \cite{Se} (for all $P\ \in
X$). Hence $\chi_P(\sigma) \in k$ is a {\em non-trivial}
$e_P^\tm$th root of unity and we obtain
\[\mbox{}-\frac{1}{e_P^\tm}\sum_{d=1}^{e_P^\tm-1} d
\cdot \chi_P^d(\sigma) = (1-\chi_P(\sigma))^{-1}\] by Lemma~(3.3)
below. Thus the character value of the right hand side at $\sigma$
is equal to
\[\sum_{P\in X^\sigma}
\frac{\Trace(\sigma|\cE(P))}{1-\chi_P(\sigma)}\] which in turn is
equal to the character value of the left hand side at $\sigma$ by
the Lefschetz fixed point formula (see Chapter~VI, \S 9 in
\cite{FL} or Example 3 in \cite{Ko1}).\\
We now fix an arbitrary subgroup $H$ of $G$ and show that, by
restricting to $H$, we obtain the corresponding formula for the
action of $H$ on $X$. This is obviously true for the left hand
sides:
\[\Res_H^G\left(\chi(G,X,\cE)\right) = \chi(H,X,\cE) \quad \textrm{in} \quad
K_0(H,k).\] In order to verify the corresponding equality for the
right hand sides, we first show the following formula for any $Q
\in Y$ and $d\in \ZZ$:
\begin{equation}
\Res_H^G \left(\Ind_{G_{\tilde{Q}}}^G\left(\cE(\tilde{Q}) \otimes
\chi_{\tilde{Q}}^d\right)\right) \cong \oplusm_{R\in X/H, \,\,
R\mapsto Q} \Ind_{H_{\tilde{R}}}^H\left(\cE(\tilde{R}) \otimes
\chi_{\tilde{R}}^d\right);
\end{equation}
here, $\tilde{Q} \in X$ is a preimage of $Q$ under $\pi$,
$\tilde{R} \in X$ is a preimage of $R$ under the canonical
projection $X \ra X/H$, $H_P$ denotes the decomposition group of
the action of $H$ on $X$ at the place $P$, and $\cE(P) \otimes
\chi_{P}^d$ denotes the obvious representation of both $G_P$ and
$H_P$ (for $P \in X$). Note that neither the left hand side nor
the right hand side of the formula~(1) depend on the chosen points
in the fibre of $\pi$ or of $X \ra X/H$. By Mackey's double coset
theorem (see Theorem~(44.2) on p.~324 in \cite{CR}), the left hand
side of (1) is equal to
\[\oplusm_{\sigma \in H\backslash G /G_{\tilde{Q}}}
\Ind_{H\cap\sigma G_{\tilde{Q}}\sigma^{-1}}^H
\left(\left(\sigma^{-1}\right)^*\left(\cE(\tilde{Q}) \otimes
\chi_{\tilde{Q}}^d\right)\right);\] furthermore, the association
$\sigma \mapsto \overline{\sigma(\tilde{Q})}$ yields a bijection
between $H\backslash G/G_{\tilde{Q}}$ and $\{R\in X/H: R\mapsto
{Q}\}$, and we have: $H\cap \sigma G_{\tilde{Q}}\sigma^{-1} = H
\cap G_{\sigma(\tilde{Q})} = H_{\sigma(\tilde{Q})}$ and
$\left(\sigma^{-1}\right)^*\Big(\cE(\tilde{Q})\otimes
\chi_{\tilde{Q}}^d\Big) \cong \cE\big(\sigma(\tilde{Q})\big)
\otimes \chi_{\sigma(\tilde{Q})}^d$. This proves formula~(1). We
now verify the above-mentioned equality for the right hand sides.
Since we already know that the $k$-dimensions agree, it suffices
to verify this equality modulo $\QQ[k[H]]$. Using formula~(1), we
obtain for the restriction of the right hand side of the formula
in Theorem~(3.1) to $H$ modulo $\QQ[k[H]]$:
\begin{eqnarray*}
\lefteqn{-\Res_H^G\left(\frac{1}{n}\sum_{P\in X}e_P^\w
\sum_{d=0}^{e_P^\tm-1} d \cdot \left[\Ind_{G_P}^G \left(\cE(P)
\otimes \chi_P^d\right)\right]\right)=}\\
&=& -\Res_H^G\left(\sum_{Q\in Y} \frac{1}{e_{\tilde{Q}}^\tm}
\sum_{d=0}^{e_{\tilde{Q}}^\tm-1} d \cdot
\left[\Ind_{G_{\tilde{Q}}}^G\left(\cE(\tilde{Q})\otimes
\chi_{\tilde{Q}}^d\right)\right]\right)\\
&=& - \sum_{Q\in Y} \frac{1}{e_{\tilde{Q}}^\tm}
\sum_{d=0}^{e_{\tilde{Q}}^\tm -1} d \sum_{R\in X/H, \,\, R\mapsto
Q} \left[\Ind_{H_{\tilde{R}}}^H\left(\cE(\tilde{R}) \otimes
\chi_{\tilde{R}}^d\right)\right]\\
&=& - \sum_{R \in X/H} \frac{1}{e_{\tilde{R}}^\tm}
\sum_{d=0}^{e_{\tilde{R}}^\tm -1} d \cdot
\left[\Ind_{H_{\tilde{R}}}^H\left(\cE(\tilde{R}) \otimes
\chi_{\tilde{R}}^d\right)\right].
\end{eqnarray*}
Let now $e_{\tilde{R}}(H)$, $e_{\tilde{R}}^\tm(H)$ and
$e_{\tilde{R}}^\w(H)$ denote the ramification indices of the cover
$X \ra X/H$ at the place $\tilde{R}$ which are defined analogously
to those of the cover $\pi$. Furthermore we put $f_{\tilde{R}} :=
e_{\tilde{R}}^\tm/e_{\tilde{R}}^\tm(H) \in \NN$. Then we can write
the latter term as follows:
\begin{eqnarray*}
\lefteqn{\mbox{}-\sum_{R\in X/H} \frac{1}{e_{\tilde{R}}^\tm}
\sum_{d=0}^{e_{\tilde{R}}^\tm(H)-1} \sum_{a=0}^{f_{\tilde{R}}}
\left(d+a \cdot e_{\tilde{R}}^\tm(H)\right) \cdot
\left[\Ind_{H_{\tilde{R}}}^H\left(\cE(\tilde{R})\otimes\chi_{\tilde{R}}^{d+a\cdot
e_{\tilde{R}}^\tm(H)}\right)\right]=}\\
&=& \mbox{}-\sum_{R\in X/H} \frac{1}{e_{\tilde{R}}^\tm}
\sum_{d=0}^{e_{\tilde{R}}^\tm(H)-1}d \cdot f_{\tilde{R}} \cdot
\left[\Ind_{H_{\tilde{R}}}^H\left(\cE(\tilde{R})\otimes\chi_{\tilde{R}}^d\right)\right]\\
&&\mbox{} - \sum_{R\in X/H} \frac{1}{e_{\tilde{R}}^\tm}
\sum_{a=0}^{f_{\tilde{R}}} a \cdot e_{\tilde{R}}^\tm(H)
\sum_{d=0}^{e_{\tilde{R}}^\tm(H)-1}
\left[\Ind_{H_{\tilde{R}}}^H\left(\cE(\tilde{R}) \otimes
\chi_{\tilde{R}}^d\right)\right].
\end{eqnarray*}
Since the first ramification group $H_{\tilde{R},1}$ is a
$p$-group (see Corollaire~3 on p.~75 in \cite{Se}) and since
$H_{\tilde{R}}$ is the semi-direct product of $H_{\tilde{R},1}$
and the cyclic group $H_{\tilde{R}}/H_{\tilde{R},1}$ (see
Corollaire~4 on p.~75 in \cite{Se}), we obtain
\[\sum_{d=0}^{e_{\tilde{R}}^\tm(H)-1}[\chi_{\tilde{R}}^d] =
[k[H_{\tilde{R}}/H_{\tilde{R},1}]] = e_{\tilde{R}}^\w(H)^{-1}
\cdot [k[H_{\tilde{R}}]]\] by Lemma~(3.4) below. Hence the latter
term is, modulo $\QQ[k[H]]$, congruent to:
\begin{eqnarray*}
\lefteqn{\mbox{} -\sum_{R\in X/H} \frac{1}{e_{\tilde{R}}^\tm(H)}
\sum_{d=0}^{e_{\tilde{R}}^\tm(H)-1} d\cdot
\left[\Ind_{H_{\tilde{R}}}^H\left(\cE(\tilde{R})\otimes\chi_{\tilde{R}}^d\right)\right]=}\\
&=& - \frac{1}{|H|} \sum_{P\in X} e_P^\w(H)
\sum_{d=0}^{e_P^\tm(H)-1} d \cdot
\left[\Ind_{H_P}^H\left(\cE(P)\otimes \chi_P^d\right)\right],
\end{eqnarray*}
as was to be shown. This completes the proof of Theorem~(3.1).
\hfill $\square$

The following easy lemma is also a crucial step in an alternative
approach to Theorem~5.2 of \cite{Ch}, see section~2a of
\cite{Er2}.

{\bf (3.3) Lemma.} {\em Let $m \in \NN$ and $\zeta \not= 1$ an
$m$th root of unity in $k$. Then we have:
\[m(\zeta-1)^{-1} = \sum_{d =1}^{m-1} d \zeta^d.\]}

{\em Proof.} $(\sum_{d=1}^{m-1} d \zeta^d)(\zeta -1) =
\sum_{d=1}^{m-1} d \zeta^{d+1} - \sum_{d =1}^{m-1} d \zeta^d =
(m-1) \zeta^m - \sum_{d=1}^{m-1} \zeta^d = m$. \hspace*{\fill}
$\square$

{\bf (3.4) Lemma.} {\em Let $H$ be the semi-direct product of a
finite $p$-group $P$ and an (arbitrary) group $C$ which acts on
$P$. Then we have in $K_0(H,k)$:
\[[k[H]] = |P| \cdot [k[C]].\]}

{\em Proof.} Let $I$ denote the augmentation ideal of the group
ring $k[P]$. The group $C$ acts on $k[P]$ in the obvious way, and
the ideals $I^r$, $r \ge 0$, of $k[P]$ are clearly $C$-stable.
Since $([\sigma]-[1])^{|P|} = [\sigma^{|P|}] - [1] = 0$ for all
$\sigma \in P$, we have $I^N = 0$ for $N$ sufficiently big.
Furthermore, the group $P$ acts trivially on the successive
quotients $I^r/I^{r+1}$, $r\ge 0$. Thus we have a finite
filtration
\[k[H] = k[P] \ast C = I^0 \ast C \supseteq I^1 \ast C \supseteq
I^2 \ast C \supseteq \ldots \supseteq I^N \ast C =0\] of the
regular representation $k[H]$ by $k[H]$-submodules such that the
successive quotients split into a direct sum of $k[H]$-modules of
the form $k[C] = k[H/P]$. (Here, $I^r \ast C$ denotes the ideal
$\oplus_{\sigma\in C}I^r[\sigma]$ in the twisted group ring
$\oplus_{\sigma \in C} k[P][\sigma] = k[P] \ast C =k[H]$.) This
proves Lemma~(3.4). \hfill $\square$

Theorem~(3.1) implies the following {\em global} relation in
$K_0(G,k)$ between the representations $\Ind_{G_P}^G(\chi_P^d)$,
$P \in X$, $d=0, \ldots, e_P^\tm -1$.

{\bf (3.5) Corollary.} {\em Let $p$ be odd. For $P \in X$ and $d
\in \{0,\ldots, e_P^\tm -1\}$ we put
\[n_{P,d}:= e_P^\w \left(d + \frac{(e_P^\w -1)(e_P^\tm+1) +
\sum_{s\ge 2}(|G_{P,s}|-1)}{2}\right) \in \NN.\] Then the element
\[\sum_{P\in X} \sum_{d=0}^{e_P^\tm -1} n_{P,d} \cdot
[\Ind_{G_P}^G(\chi_P^d)] \in K_0(G,k)\] is divisible by $n = |G|$
in $K_0(G,k)$.}

{\em Proof.} Apply Theorem~(3.1) to $\cE= \cO_X$ and use the
equality
\[[k[G]] = e_P^\w \sum_{d=0}^{e_P^\tm-1} \left[\Ind_{G_P}^G(\chi_P^d)\right]\]
which follows from Lemma~(3.4). \hfill $\square$

\bigskip

\section*{\S 4 Galois Structure in the Weakly Ramified Case}

In this section we will generalize several results of Kani and
Nakajima on the Galois module structure of Zariski cohomology
groups of curves from the tamely ramified to the weakly ramified
case.

We keep the assumptions and notations introduced in \S 2 and \S 3.
In addition we assume in this section that the cover $\pi: X \ra
Y$ is weakly ramified.

We begin with recalling the following crucial properties of weakly
ramified covers.

{\bf (4.1) Lemma.} {\em For any $P \in X$, the first ramification
group $G_{P,1}$ is an abelian group of exponent $p$, the factor
group $G_P/G_{P,1}$ is cyclic of order prime to $p$ and the
natural action of $G_P/G_{P,1}$ on $G_{P,1} \backslash \{1\}$ is
free. In particular, $G_P$ is the semidirect product of $G_{P,1}$
and $G_P/G_{P,1}$ and we have: $e_P^\w \equiv 1$ mod $e_P^\tm$. }

{\em Proof.} This is proved on the pages 74-77 in Serre's book
\cite{Se}, see in particular Proposition~9 and the corollaries of
Proposition~7. \hfill $\square$

{\bf (4.2) Lemma.} {\em Let $H$ be the semi-direct product of a
finite $p$-group $P$ with a finite group $C$ which acts an $P$. We
assume that the action of $C$ on $P \backslash \{1\}$ is free, and
we put $a:=(|P|-1)/|C| \in \NN$. Furthermore let $V$ be a
$k$-representation of $C$ (of finite dimension) which we view also
as a $k$-representation of $H$ via the canonical projection $H \ra
C$. Then we have:\\
(a) The induced representation $\Ind_C^H(V)$ is the
$k[H]$-projective cover of $V$.\\
(b) $[\Ind_C^H(V)] = [V] + a \cdot \dim_k(V) \cdot [k[C]]$ in
$K_0(H,k)$. }

{\em Proof.} The order of $C$ is prime to $p$ by assumption. Hence
$V$ is a projective $k[C]$-module and $\Ind_C^H(V)$ is a
projective $k[H]$-module. Furthermore we have an obvious
$k[H]$-epimorphism $\Ind_C^H(V) \ra V$. Since $V \cong
\Ind_C^H(V)/{\rm rad}(k[H])\Ind_C^H(V)$, the $k[H]$-module
$\Ind_C^H(V)$ is minimal with these properties. Thus we have
proved part~(a). We now prove part~(b). For all $x\in P$ and
$\xi,\eta \in C$ we obviously have:
\[(x,\xi)\cdot (1,\eta) \cdot (x,\xi)^{-1} = (x \cdot (\xi \eta
\xi^{-1})(x^{-1}), \xi \eta \xi^{-1}) \quad \textrm{in} \quad
P\rtimes C = H.\] Hence, by assumption on the action of $C$ on
$P$, the intersection $C\cap \sigma C \sigma^{-1}$ is trivial for
all $\sigma \in H \backslash C$. Using Mackey's double coset
theorem (see Theorem~(44.2) on p.~324 in \cite{CR}), we thus
obtain:
\begin{eqnarray*}
\lefteqn{\Res_C^H\left(\Ind_C^H(V)\right) \cong \oplusm_{\sigma
\in C\backslash G/C} \Ind_{C\cap \sigma C \sigma^{-1}}^C
\left(\Res_{C\cap \sigma C
\sigma^{-1}}^C\left((\sigma^{-1})^*(V)\right)\right)} \hspace*{20em}\\
&=& V \oplus \left(\oplusm^{a \cdot \dim_k(V)} k[C] \right).
\end{eqnarray*}
Furthermore the restriction homomorphism $\Res_C^H : K_0(H,k) \ra
K_0(C,k)$ is bijective since the irreducible $k$-representations
of $C$ considered as $k$-representations of $H$ are the only
irreducible $k$-representations of $H$ (because $k[H]/{\rm
rad}(k[H]) \cong k[C]$). This proves part~(b). \hfill $\square$

In the sequel, the $k[G_P]$-projective cover of the
$k[G_P]$-module $\chi_P^d$ is denoted by $\Cov(\chi_P^d)$ (for all
$P \in X$ and $d \in \ZZ$). The following theorem gives a {\em
global} relation between the projective $k[G]$-modules
$\Ind_{G_P}^G\left(\Cov(\chi_P^d)\right)$, $P\in X$, $d=1, \ldots,
e_P^\tm-1$.

{\bf (4.3) Theorem.} {\em There is a (unique) projective
$k[G]$-module $N$ such that
\[\oplusm^n N \cong \oplusm_{P\in X} \oplusm_{d=1}^{e_P^\tm-1}
\oplusm^{e_P^\w \cdot d}
\Ind_{G_P}^G\left(\Cov(\chi_P^d)\right).\]}

{\em (4.4) Remark.} If $\pi$ is assumed to be not only weakly but
tamely ramified, then the projective $k[G]$-module on the right
hand side has the following somewhat simpler shape:
\[\oplusm_{P\in X} \oplusm_{d=1}^{e_P-1} \oplusm^d
\Ind_{G_P}^G(\chi_P^d).\] In particular, Theorem~(4.3) generalizes
the first part of Theorem~2 in Kani's paper~\cite{Ka} and
Theorem~2(i) in Nakajima's paper~\cite{Na2}.

{\em Proof of Theorem~(4.3).} Let $E$ denote the equivariant
divisor $E := \sum_{P\in X}(e_P^\w-1)[P]$ on $X$. Then we have in
$K_0(G,k)_\QQ$:
\begin{eqnarray*}
\lefteqn{\chi(G,X,\cO_X(E)) = }\\
&=& (1-g_Y)[k[G]]\\
&&\mbox{} + \frac{1}{n} \sum_{P\in X} \left((e_P^\w-1) -
\frac{(e_P^\w-1)(e_P^\tm+1)}{2} \right) [k[G]]\\
&& \mbox{} - \frac{1}{n} \sum_{P\in X} e_P^\w
\sum_{d=1}^{e_P^\tm-1} d \cdot
\left[\Ind_{G_P}^G\left(\chi_P^{d+1-e_P^\w}\right)\right]
\hspace*{7em}
\textrm{(by Theorem~(3.1))}\\
&=& (1-g_Y)[k[G]]\\
&& \mbox{} - \frac{1}{n} \sum_{P\in X}
\frac{(e_P^\tm-1)(e_P^\w-1)}{2} [k[G]] \\
&& \mbox{} - \frac{1}{n} \sum_{P\in X} e_P^\w
\sum_{d=0}^{e_P^\tm-1} d  \cdot
\left[\Ind_{G_P}^G(\chi_P^d)\right] \hspace*{2em} \textrm{(since
}e_P^\w \equiv 1 \textrm{ mod } e_P^\tm \textrm{ by
Lemma~(4.1))}\\
&=& (1-g_Y) [k[G]] \\
&& \mbox{} - \frac{1}{n} \sum_{P\in X}
\frac{(e_P^\tm-1)e_P^\tm}{2} \cdot \frac{e_P^\w-1}{e_P^\tm} [k[G]]\\
&& \mbox{} - \frac{1}{n} \sum_{P\in X} e_P^\w
\sum_{d=0}^{e_P^\tm-1} d \cdot
\left(\left[\Ind_{G_P}^G\left(\Cov(\chi_P^d)\right)\right] -
\frac{e_P^\w-1}{e_P^\tm}\left[\Ind_{G_P}^G\left(k[G_P/G_{P,1}]\right)\right]\right)\\
&&\hspace*{18em} \textrm{(by Lemma~(4.1) and Lemma~(4.2))}\\
&=& (1-g_Y) [k[G]]\\
&& \mbox{} - \frac{1}{n} \sum_{P\in X} e_P^\w
\sum_{d=0}^{e_P^\tm-1} d \cdot
\left[\Ind_{G_P}^G\left(\Cov(\chi_P^d)\right)\right] \hspace*{1em}
\textrm{(by Lemma~(4.1) and Lemma~(3.4))}.
\end{eqnarray*}
Hence, by Theorem~(2.1)(a)(i), the class of the projective
$k[G]$-module
\[\oplusm_{P\in X} \oplusm_{d=1}^{e_P^\tm-1} \oplusm^{d \cdot
e_P^\w} \Ind_{G_P}^G\left(\Cov(\chi_P^d)\right)\] in $K_0(k[G])$
is divisible by $n$ in $K_0(k[G])$. Writing the elements of
$K_0(k[G])$ as integral linear combinations of a basis of
$K_0(k[G])$ consisting of indecomposable projective
$k[G]$-modules, we see that the quotient is again the class of a
projective $k[G]$-module, say $N$. Since two projective
$k[G]$-modules, whose classes in $K_0(k[G])$ are equal, are
already isomorphic, Theorem~(4.3) is now proved. \hfill $\square$

The following theorem expresses the equivariant Euler
characteristic $\chi(G,X,\cO_X(D))$ as an  {\em integral} linear
combination of classes of explicit projective $k[G]$-modules (for
any equivariant divisor $D$ as in Theorem(2.1)(a)).

{\bf (4.5) Theorem.} {\em Let $D=\sum_{P\in X} n_P[P]$ be an
equivariant divisor on $X$ with $n_P \equiv -1$ mod $e_P^\w$ for
all $P \in X$. For any $P \in X$, we write
\[n_P = (e_P^\w-1) + (l_P +m_P e_P^\tm) e_P^\w\]
with $l_P \in \{0, \ldots, e_P^\tm-1\}$ and $m_P \in \ZZ$.
Furthermore, for any $Q\in Y$, we choose a preimage $\tilde{Q}\in
X$ of $Q$ under $\pi$. Then we have in $K_0(k[G])$:
\[\chi(G,X,\cO_X(D)) = - [N] + \sum_{Q\in Y}
\sum_{d=1}^{l_{\tilde{Q}}}
\left[\Ind_{G_{\tilde{Q}}}^G\left(\Cov(\chi_{\tilde{Q}}^{-d})\right)\right]
+ \left(1-g_Y + \sum_{Q\in Y} m_{\tilde{Q}}\right)[k[G]].\]}

{\em (4.6) Remark.} Note that $l_P = 0 = m_P$ for all but finitely
many $P\in X$. If $\pi$ is tamely ramified at $P$, then $l_P$ is
obviously the unique number in $\{0, \ldots, e_P^\tm-1\}$ such
that $\cO_X(D)(P) \cong \chi_P^{-l_P}$ (as $k[G_P]$-modules). If
$\pi$ is tamely ramified everywhere, then Theorem~(4.5) implies
the congruence
\[\chi(G,X,\cO_X(D)) \equiv -[N] + \sum_{Q\in Y}
\sum_{d=1}^{l_{\tilde{Q}}}
\left[\Ind_{G_{\tilde{Q}}}^G(\chi_{\tilde{Q}}^{-d})\right]
\textrm{ mod } \ZZ[k[G]]\] for an arbitrary equivariant divisor
$D$ on $X$. This congruence is a reformulation of Theorem~2(ii) in
\cite{Na2} applied to $\cE:= \cO_X(D)$.

{\em Proof of Theorem~(4.5).} Let $E$ denote the divisor
$\sum_{P\in X} (e_P^\w-1)[P]$ as in the proof of Theorem~(4.3). We
first compute the difference $\chi(G,X,\cO_X(D)) -
\chi(G,X,\cO_X(E))$ in $K_0(G,k)_\QQ$.
\begin{eqnarray*}
\lefteqn{\chi\left(G,X,\cO_X(D)\right) -
\chi\left(G,X,\cO_X(E)\right) = } \\
&=& \frac{1}{n} \sum_{P\in X} n_P [k[G]] - \frac{1}{n} \sum_{P\in
X}(e_P^\w-1)[k[G]] \\
&& \mbox{} - \frac{1}{n} \sum_{P\in X} e_P^\w
\sum_{d=1}^{e_P^\tm-1} d \,
\left[\Ind_{G_P}^G\left(\chi_P^{d-n_P}\right)\right] + \frac{1}{n}
\sum_{P\in X} e_P^\w \sum_{d=1}^{e_P^\tm-1}
d \, \left[\Ind_{G_P}^G\left(\chi_P^{d+1-e_P^\w}\right)\right]\\
&& \hspace*{25em} \textrm{(by Theorem~(3.1))}\\
&=& \frac{1}{n} \sum_{P\in X}(l_P+m_Pe_P^\tm)e_P^\w[k[G]]\\
&& \mbox{} - \frac{1}{n} \sum_{P\in X}
e_P^\w\left(\sum_{d=0}^{e_P^\tm-1} d \,
\left[\Ind_{G_P}^G\left(\chi_P^{d-l_P}\right)\right]
-\sum_{d=0}^{e_P^\tm-1} d \,
\left[\Ind_{G_P}^G\left(\chi_P^d\right)\right]\right)\\
&&\hspace*{17em} (\textrm{since } e_P^\w \equiv 1 \textrm{ mod }
e_P^\tm \textrm{ by Lemma~(4.1)})\\
&=& \frac{1}{n}\sum_{P\in X} m_Pe_P[k[G]] + \frac{1}{n} \sum_{P\in
X} l_P [k[G]] + \frac{1}{n} \sum_{P\in X} l_P(e_P^\w-1) [k[G]]\\
&& \mbox{} - \frac{1}{n} \sum_{P\in X} e_P^\w
\left(\sum_{d=0}^{e_P^\tm-1}
l_P\left[\Ind_{G_P}^G(\chi_P^d)\right]-\sum_{d=e_P^\tm-l_P}^{e_P^\tm-1}
e_P^\tm \left[\Ind_{G_P}^G(\chi_P^d)\right]\right)\\
&=& \sum_{P\in X} \frac{e_P}{n} m_P [k[G]] + \sum_{P\in X}
\frac{e_P}{n} \, l_P \, \frac{e_P^\w-1}{e_P^\tm}
\left[\Ind_{G_P}^G(k[G_P/G_{P,1}])\right] \\
&& \mbox{} + \sum_{P\in X} \frac{e_P}{n}
\sum_{d=e_P^\tm-l_P}^{e_P^\tm-1}\left[\Ind_{G_P}^G(\chi_P^d)\right]
\hspace*{2em}\\
&& \hspace*{14em} (\textrm{since }[k[G_P]] = e_P^\w
\sum_{d=0}^{e_P^\tm-1} [\chi_P^d] \textrm{ by Lemma~(3.4)})\\
&=& \sum_{Q\in Y} m_{\tilde{Q}}[k[G]] + \sum_{Q\in
Y}\sum_{d=1}^{l_{\tilde{Q}}}
\left[\Ind_{G_{\tilde{Q}}}^G\left(\Cov(\chi_{\tilde{Q}}^{-d})\right)\right]\\
&&\hspace*{18em}
 \textrm{ (by Lemma~(4.1) and Lemma~(4.2))}.
\end{eqnarray*}
This result for the difference $\chi(G,X,\cO_X(D)) -
\chi(G,X,\cO(E))$ together with the formula
\[\chi(G,X,\cO_X(E)) = (1-g_Y)[k[G]] - [N]\]
from the proof of Theorem~(4.3) now obviously implies the desired
formula in Theorem~(4.5). \hfill $\square$

The ideal sheaf of the reduced effective divisor on $X$ consisting
of all (wildly) ramified points plays a central role in Pink's
paper \cite{Pi}. The following corollary computes the Galois
module structure of its first cohomology group.

{\bf (4.7) Corollary.} {\em Let $S$ be a $G$-stable non-empty
finite set of closed points on $X$ which contains all {\em wildly
ramified} points, i.e.\ all points $P \in X$ with $e_P^\w \not=
1$. Let $\cI(S)$ denote the ideal sheaf of $S$. Then the
$k[G]$-module $H^1(X,\cI(S))$ is stably isomorphic to $N \oplus
\left(\oplus_{Q\in \pi(S)}
\Ind_{G_{\tilde{Q}}}^G\left(\Cov(\chi_{\tilde{Q}}^0)\right)\right)$.
More precisely we have:
\[H^1(X,\cI(S)) \oplus k[G] \cong N \oplus \left(\oplusm_{Q\in
\pi(S)}
\Ind_{G_{\tilde{Q}}}^G\left(\Cov(\chi_{\tilde{Q}}^0)\right)\right)
\oplus \left(\oplusm^{g_Y} k[G]\right).\]}

{\em Proof.} Since $S \not= \emptyset$, we have $H^0(X,\cI(S)) =
0$. Hence we obtain the following equality in $K_0(G,k)$:
\begin{eqnarray*}
\lefteqn{\left[H^1(X,\cI(S))\right] = -
\chi \left(G,X,\cO_X\left(\sum_{P\in S} - [P]\right)\right) =}\\
&=& [N] - \sum_{Q\in\pi(S)} \sum_{d=1}^{e_{\tilde{Q}}^\tm-1}
\left[\Ind_{G_{\tilde{Q}}}^G\left(\Cov(\chi_{\tilde{Q}}^{-d})\right)\right]
+ \sum_{Q\in \pi(S)} [k[G]] + (g_Y-1) [k[G]]\\
&& \hspace*{25em} \textrm{(by Theorem~(4.5))}\\
&=& [N] + \sum_{Q\in \pi(S)}
\left[\Ind_{G_{\tilde{Q}}}^G\left(\Cov(\chi_{\tilde{Q}}^0)\right)\right]
+ (g_Y-1) [k[G]]\\
&& \hspace*{18em} \textrm{(by Lemma~(4.1) and Lemma~(4.2)).}
\end{eqnarray*}
Furthermore, $H^1(X,\cI(S))$ is a projective $k[G]$-module by
Theorem~(2.1)(a)(ii). Since two projective $k[G]$-modules, whose
classes are equal in $K_0(G,k)$, are already isomorphic,
Corollary~(4.7) is now proved. \hfill $\square$

The following corollary computes the Galois module structure of
the space of global meromorphic differentials on $X$ which are
logarithmic along all ramified points. It generalizes the second
part of Theorem~2 in Kani's paper \cite{Ka} from the tamely
ramified to the weakly ramified case.

{\bf (4.8) Corollary.} {\em Let $S$ be a $G$-stable non-empty
finite set of closed points on $X$ which contains all ramified
points. Let $\Omega_X(S)$ denote the sheaf of meromorphic
differentials on $X$ which are logarithmic along $S$. Then the
$k[G]$-module $H^0(X,\Omega_X(S)) \oplus N$ is free of rank
$|S/G|+g_Y-1$.}

{\em Proof.} We use the notations introduced in the proof of
Corollary~(2.3). Then the divisor $D= K_X + \sum_{P\in S} [P]$
satisfies the condition of Theorem~(4.5) and the corresponding
integers $l_P$, $m_P$, $P \in X$, are given as follows: $l_P =0$
for all $P \in X$, $m_P = r_{\pi(P)} +1$ for $P \in S$ and $m_P =
r_{\pi(P)}$ for $P \in X \backslash S$. Thus we obtain:
\begin{eqnarray*}
\lefteqn{[H^0(X,\Omega_X(S)) \oplus N] = [H^0(X,\cO_X(D))] + [N]
=}\\
&=& (1-g_Y + \sum_{Q \in Y} r_Q + |S/G|) [k[G]] \quad \textrm{(by
Theorem~(4.5))}\\
&=& (g_Y-1 + |S/G|) [k[G]] \quad \textrm{in} \quad K_0(k[G]) \quad
\textrm{(since } \deg(\Omega_Y) = 2g_Y-2).
\end{eqnarray*}
As in Corollary~(4.7), this implies Corollary~(4.8). \hfill
$\square$

{\em (4.9) Remark.} Alternatively, Corollary~(4.8) can be derived
from Corollary~(4.7) using the Serre duality isomorphism
$H^0(X,\Omega_X(S)) \cong H^1(X,\cI(S))^*$ and the isomorphism
\[N^* \oplus \left( \oplusm_{Q\in \pi(S)} \Ind_{G_{\tilde{Q}}}^G
\left(\Cov(\chi_{\tilde{Q}}^0)\right)^*\right) \oplus N \cong
\oplusm^{|S/G|} k[G]\] which may easily be checked.

\bigskip

Faculty of Mathematical Studies\\
University of Southampton\\
Southampton SO17 1BJ\\
United Kingdom\\
{\em e-mail:} bk@maths.soton.ac.uk.

\end{document}